\documentclass[
12pts ] {article}

\usepackage[english]{babel}
\usepackage{graphicx,color}
\usepackage{amsfonts}
\usepackage{amssymb,amsmath,amsthm,scalefnt,mathrsfs}
\usepackage[all]{xy}

\title{Right-Permutative Cellular Automata on Topological Markov Chains}

\author{Marcelo Sobottka \footnote{Departamento de Ingenier\'{\i}a
Matem\'{a}tica, Facultad de Ciencias F\'{\i}sicas y Matem\'{a}ticas, Universidad
de Concepci\'{o}n, casilla 160-C, correo 3, Concepci\'{o}n, Chile (e-mail:
msobottka@udec.cl). This work was supported by {\sc MECESUP UCH0009}
and {\sc Nucleus Millennium Information and Randomness ICM
P04-069-F}. Part of this work was carried out while the author was
invited by the research group {\sc Modelos Estoc\'{a}sticos (DGA)} of
the Universidad de Zaragoza.}}
\date{\ }
\date{\ }

\newtheorem{theo}{Theorem}[section]
\newtheorem{cor}[theo]{Corollary}
\newtheorem{lem}[theo]{Lemma}

\newtheorem{prop}[theo]{Proposition}
\newtheorem{defn}[theo]{Definition}
\newtheorem{rem}[theo]{Remark}
\newtheorem{ex}[theo]{Example}
\numberwithin{equation}{section}

\newcommand{\abs}[1]{\left\vert#1\right\vert}
\newcommand{\set}[1]{\left\{#1\right\}}

\newcommand{\Z}{\mathbb{Z}}
\newcommand{\N}{\mathbb{N}}

\newcommand{\s}{\sigma}

\newcommand{\F}{\mathcal{F}}

\newcommand{\Pp}{\mathcal{P}}

\newcommand{\n}{^{-1}}

\newcommand{\cqd}{\begin{flushright} $\square$\\\end{flushright}}

\newcommand{\Sgg}{\mathfrak{G}}
\newcommand{\h}{\mathfrak{H}}
\newcommand{\GZ}{G^\mathbb{Z}}
\newcommand{\e}{\mathbf{e}}

\begin{document}
\maketitle

\setcounter{theo}{0}


\begin{abstract}

In this paper we consider cellular automata $(\Sgg,\Phi)$ with
algebraic local rules and such that $\Sgg$ is a topological Markov
chain which has a structure compatible to this local rule. We
characterize such cellular automata and study the convergence of the
Ces\`{a}ro mean distribution of the iterates of any probability measure
with complete connections and summable decay.\\

{\bf Subj-class:} DS, PR, CO.

{\bf MSC-class:} Primary: 37B15, 54H20; Secondary: 28D99, 37B20.

\end{abstract}

\bigskip
\hrule
\noindent
{\footnotesize\em This is a pre-copy-editing, author-produced PDF of an article accepted for publication in Discrete and Continuous Dynamical Systems - Series A (DCDS-A), following peer review. The definitive publisher-authenticated version Marcelo Sobottka, Topological quasi-group shifts. Disc. and Cont. Dynamic. Systems (2008), 20, 4, 1095-1109, is available online at: http://www.aimsciences.org/journals/displayArticles.jsp?paperID=3147 .}
\hrule
\bigskip

\section{Introduction}\label{1}

Let $\GZ$ be the two sided full shift on the finite alphabet $G$,
and $\s:\GZ\to\GZ$ be the shift map. Suppose $\Sgg\subseteq\GZ$ is a
topological Markov Chain which, without loss of generality, we can
consider uses all alphabet $G$.

Consider the cellular automaton $(\Sgg,\Phi)$ which has a local rule
defined from some algebraic operation on $G$. Motivated by their
several applications in information theory, physics, and biological
sciences, among others, the problem of to characterize and to
analyze the dynamical behavior of such cellular automata has been
widely investigated. More specifically, there are three important
questions about $(\Sgg,\Phi)$: if it is possible to recode it in the
way to understand and to classify its dynamics (see
\cite{HostMaassMartinez} and \cite{mms}); what $\s$-invariant
probability measures are also $\Phi$-invariant (see
\cite{HostMaassMartinez}, \cite{Pivato} and \cite{Sablik}); and how
$\s$-invariant probability measures evolve under the dynamics
of $\Phi$ (see \cite{HostMaassMartinez}, \cite{mmpy} and \cite{mms}).\\

When $\Sgg=\GZ$ and $(\Sgg,\Phi)$ is a right-permutative
$\Psi$-associative or $N$-scaling cellular automaton,
Host-Maass-Mart\'{\i}nez \cite{HostMaassMartinez} proved that it is
topologically conjugate to an affine cellular automaton product a
translation $(K^\Z\times B^\Z,\Phi_K\times\mathbf{G}_B)$. Moreover,
they showed sufficient conditions under which the unique
shift-affine invariant measure is the maximum entropy measure
(property which is known as rigidity), and studied the convergence
of the Ces\`{a}ro mean distribution of $\s$-invariant probability
measures under the action of $\Phi$. The results of
\cite{HostMaassMartinez} about rigidity were generalized by Pivato
\cite{Pivato} for the case of bipermutative endomorphic cellular
automata. Is his work, Pivato also showed results about the
characterization of the topological dynamics of bipermutative
cellular automata. Later, the rigidity results of Pivato were
generalized by Sablik \cite{Sablik} who also includes the case of
$\Sgg$ being a proper subgroup shift of $\GZ$.

Recently, Mass-Mart\'{\i}nez-Sobottka \cite{mms} have showed that if
$(\Sgg,+)$ is an Abelian subgroup shift and a $p^s$-torsion for some
prime number $p$, and $\Phi$ is an affine cellular automaton given
by $\Phi:=\mathbf{a}\cdot id\ +\ \mathbf{b}\cdot\s\ +\ \mathbf{c}$,
where $\mathbf{a},\mathbf{b}\in\N$ are relatively prime to $p$, and
$\mathbf{c}\in\Sgg$ is a constant sequence, then the Ces\`{a}ro mean
distribution of any measure with complete connections (compatible
with $\Sgg$) and summable decay under the action of $\Phi$ converges
to the maximum entropy measure on $\Sgg$. The proof of this result
combines regeneration theory, combinatorics, and the recodification
of $\Sgg$. As consequence of the convergence of the Ces\`{a}ro mean
distribution we get a rigidity property, namely: the unique
$(\s,\Phi)$-invariant measure with complete connections (compatible
with $\Sgg$) and summable decay
for that case is the maximum entropy measure.\\

This paper concentrates mainly on the first problem, characterizing
the dynamical behavior of bipermutative and some right-permutative
cellular automata defined on subshifts $\Sgg$ which are not
necessarily subgroup shifts, but which have some algebraic
structure. As a direct application of these results we recuperate
several results about rigidity and about the evolution of
$\s$-invariant measures under the action of $\Phi$.

This paper is organized as follows. In {\S}\ref{2} we develop the
background. In {\S}\ref{3}, we define the class of
structurally-compatible cellular automata and study the case of
bipermutative cellular automata. In {\S}\ref{4} we study the
representation of right-permutative $\Psi$-associative or
$N$-scaling cellular automata. In {\S}\ref{5} we present some
sufficient conditions under which a block code preserves the
properties of complete connections and summable decay of a
probability measure, and so in {\S}\ref{6} we apply the results
obtained in the previous sections to study the convergence of the
Ces\`{a}ro mean distribution. In {\S}\ref{7} we gives some results about
rigidity.

\section{Background}\label{2}

Let $\Sgg\subseteq\GZ$ be a subshift. Given $\mathbf{g}\in\Sgg$, and
$m\leq n$, we denote by $\mathbf{g}_m^n=(g_m,g_{m+1},\ldots,g_n)$.
For $k\geq 1$, denote by $\Sgg_k$ the set of all allowed words with
length $k$ in $\Sgg$. Given $\mathbf{g}\in \Sgg_k$,
$\mathbf{g}=(g_1\ldots,g_k)$ we write $\F(\mathbf{g})$, as the {\em
follower set} of $\mathbf{g}$ in $\Sgg$:
$$\F(\mathbf{g})=\set{h\in G:(g_1,\ldots,g_k,h)\in\Sgg_{k+1}}.$$

In the same way, we define $\Pp(\mathbf{g})$ the set of predecessors
of $\mathbf{g}\in\Sgg_k$ in $\Sgg$.

We say a subshift $\Sgg$ is a {\em topological Markov chain} if for
any $k\geq 1$ and $\mathbf{g}=(g_1,\ldots,g_k)\in\Sgg_k$ we have
$\F(\mathbf{g})=\F(g_k)$, which means $\Sgg$ can be thought as
generated by a bi-infinite walking on an oriented graph. A
topological Markov chain $\Sgg$ is {\em irreducible} if and only if
for any $u,w\in G$ there exist $k\geq 1$ and
$(v_1,\ldots,v_k)\in\Sgg_k$ such that
$(u,v_1,\ldots,v_k,w)\in\Sgg_{k+2}$, and it is {\em mixing} if there
exists $q\geq 1$ such that for any $k\geq q$ and $u,w\in G$ we
always can find $(v_1,\ldots,v_k)\in\Sgg_k$ such that
$(u,v_1,\ldots,v_k,w)\in\Sgg_{k+2}$.\\

Denote by $\Sgg^-$ and $\Sgg ^+$ , the projections of $\Sgg$ on
$G^{-\N^*}$ and $G^{\N}$ respectively. Given $w\in\Sgg^-$ denote by
$\Sgg_w^+$ the projection on $\Sgg^+$ of the set of all
sequences $(g_i)_{i\in\Z}\in\Sgg$, with $g_i=w_i$ for $i\leq -1$.\\

Let $\s:\Sgg\to\Sgg$ be the {\em shift map}, which is defined for
every $\mathbf{g}\in\Sgg$ and $n\in\Z$ as
$(\s(\mathbf{g}))_n=g_{n+1}$.

We say a map $\Theta:\Lambda\to\Lambda'$, between two topological
Markov chains is a $(\ell+r+1)$-block code if it has a local rule
$\theta:\Lambda_{\ell+r+1}\to \Lambda'_1$ such that for any
$\mathbf{x}=(x_i)_{i\in\Z}\in\Lambda$ and $j\in\Z$ follows that
$\bigl(\Theta(\mathbf{x})\bigr)_j=\theta(x_{j-\ell},\ldots,x_{j+r})$.
Under these notations, we say $\Theta$ has memory $\ell$ and
anticipation $r$. We recall a map $\Theta$ is a block code if and
only if it is continuous and commutes with the shift map.

$$\xymatrix{
\mathbf{x}       & = & (  \ldots  , & *+[F]{x_{j-\ell},\ldots,x_j\ldots,x_{j+r}}\ar[d]^{\theta} & ,  \ldots  ) \\
\Theta(\mathbf{x}) & = & (  \ldots  , &
\bigl(\Theta(\mathbf{x})\bigr)_j & , \ldots  ) }$$\\ 

A {\em cellular automaton (c.a.)} is a pair $(\Sgg,\Phi)$, where
$\Phi:\Sgg\to\Sgg$ is $(\ell+r+1)$-block code. Without loss of
generality we always can consider $\ell=0$ and so to say the c.a.
has radio $r$.

A c.a. with radio $r$ is said {\em right permutative}, if its local
rule $\phi$ verifies for any fixed word
$(w_0,\ldots,w_{r-1})\in\Sgg_r$ that the map $g\mapsto
\phi(w_0,\ldots,w_{r-1},g)$ is a permutation on $G$. In the
analogous way we define left permutativity. When a c.a. is right and
left permutative, we will say it is {\em bipermutative}. From now
on, we will consider that $\Phi:\Sgg\to\Sgg$ is a restriction on
$\Sgg$ of some c.a. $\tilde{\Phi}:\GZ\to\GZ$. It is equivalent to
say that there exists a map $\tilde{\phi}:G^{r+1}\to G$ such that
the local rule of $\Phi$ is $\phi=\tilde{\phi}|_{\Sgg_{r+1}}$.

Let us to define three types of cellular automata which are
fundamental in this work:
                {\bf translations:} $(\Sgg,\mathbf{g})$ is a translation if
                $\mathbf{g}=\mathbf{s}\circ\s$,
                where $\mathbf{s}:\GZ\to \GZ$ is a 1-block code with
                local rule $s:G\to G$ which is a permutation on $G$;
                {\bf affine c.a.:} $(\Sgg,\Phi)$ is an affine c.a. if
                its local rule is given by $\phi(a,b)=\eta(a)+\rho(b)+c$, where
                $+$ is an Abelian group operation on $G$, $\eta:G\to G$
                and $\rho:G\to G$ are two commuting automorphisms (that is,
                $\eta\circ\rho=\rho\circ\eta$), and $c\in G$; and
                {\bf group c.a.:} $(\Sgg,\Phi)$ is a group c.a. if
                its local rule is given by $\phi(a,b)= a+b$, where
                $+$ is an Abelian group operation on $G$.\\

We say a binary operation $*$ on $\Sgg$ is $(\ell+r+1)$-block if the
map $(\mathbf{x},\mathbf{y})\in\Sgg\times\Sgg\mapsto
\mathbf{x}*\mathbf{y}\in\Sgg$ is a $(\ell+r+1)$-block code. When $*$
is a (quasi) group operation, then we say $(\Sgg,*)$ is a (quasi)
group
shift.\\

Let $\mu$ be any $\s$-invariant probability measure on $\Sgg$. For a
past $w\in \Sgg^-$, $w=(\ldots,w_{-2},w_{-1})$, let $\mu_w$ be the
probability measure on $\Sgg_w^+$ obtained for $\mu$ conditioning to
the past $w$.

We say $\mu$ has {\em complete connections} (compatible with $\Sgg$)
if given $a\in G$, for all $w\in \Sgg^-$ such that $a\in
\F(w_{-1})$, one has $\mu_w(a)>0$.

If $\mu$ is a probability measure with complete connections, we
define the quantities $\gamma_m$, for $m\geq 1$, by

$$\gamma_m:=\sup\set{\abs{\frac{\mu_v(a)}{\mu_w(a)}-1}:\quad \begin{array}{l}v,w\in
\Sgg^-;\quad v_{-i}=w_{-i},\ 1\leq i\leq m;\\
a\in\F(v_{-1})=\F(w_{-1})\end{array}}.$$\\

When $\sum_{m\geq 1}\gamma_m< \infty$, we say $\mu$ has {\em
summable decay}.\\ 



\section{Cellular automata with algebraic local rules}\label{3}

In this section we shall define the class of {\sc
structurally-compatible} cellular automata, which is the subject of
this work. Moreover, we will study the case of
structurally-compatible
bipermutative c.a..\\

{\defn We say a cellular automaton $(\Sgg,\Phi)$ with radio $1$ is
structurally compatible (SC) if it verifies the following property:

\begin{equation}\label{IV2}
(x_i)_{i\in\Z},(y_i)_{i\in\Z}\in \Sgg\Longrightarrow
\bigl(\phi(x_i,y_i)\bigr)_{i\in\Z}\in\Sgg,
\end{equation}

where $\phi$ denotes the local rule of $\Phi$.}\\

Define $\bullet$ as the binary operation on $G$ giving for all
$a,b\in G$ by $a\bullet b:=\phi(a,b)$. The structural compatibility
implies we can consider the componentwise operation $*$ on $\Sgg$:

$$\forall
(x_i)_{i\in\Z},(y_i)_{i\in\Z}\in\Sgg,(x_i)_{i\in\Z}*(y_i)_{i\in\Z}:=(x_i\bullet
y_i)_{i\in\Z},$$\\

Notice neither $\bullet$ nor $*$ are necessarily algebraic
operations on $G$ and $\Sgg$ respectively. However, the c.a. is left
permutative (as well right permutative or bipermutative) if and only
if $\bullet$ (and so $*$) is a left cancellable operation (as well
right cancellable or a quasi-group operation respectively). We
recall that an operation which is left-right cancellable is called a
quasi-group operation.

In terms of $*$, the map $\Phi$ can be written  as
$$\Phi=id*\s.$$\\

{\ex\label{IV3} Let $\bullet$ be the quasi-group operation
on $G=\set{a_i,b_i,c_i,d_i:\ i=1,2,3}$, giving by the following Latin square:\\

\begin{center}
\begin{tabular}{|c||c c c c c c c c c c c c|}
\hline

$\bullet$

      & $a_1$ & $b_1$ & $c_1$ & $d_1$ & $a_2$ & $b_2$ & $c_2$ & $d_2$ & $a_3$ & $b_3$ & $c_3$ & $d_3$ \\
\hline\hline

$a_1$ & $a_3$ & $b_3$ & $c_3$ & $d_3$ & $a_2$ & $b_2$ & $c_2$ & $d_2$ & $a_1$ & $b_1$ & $c_1$ & $d_1$ \\

$b_1$ & $b_3$ & $a_3$ & $d_3$ & $c_3$ & $b_2$ & $a_2$ & $d_2$ & $c_2$ & $b_1$ & $a_1$ & $d_1$ & $c_1$ \\

$c_1$ & $c_3$ & $d_3$ & $a_3$ & $b_3$ & $c_2$ & $d_2$ & $a_2$ & $b_2$ & $c_1$ & $d_1$ & $a_1$ & $b_1$ \\

$d_1$ & $d_3$ & $c_3$ & $b_3$ & $a_3$ & $d_2$ & $c_2$ & $b_2$ & $a_2$ & $d_1$ & $c_1$ & $b_1$ & $a_1$ \\

$a_2$ & $a_2$ & $b_2$ & $c_2$ & $d_2$ & $a_1$ & $b_1$ & $c_1$ & $d_1$ & $a_3$ & $b_3$ & $c_3$ & $d_3$ \\

$b_2$ & $b_2$ & $a_2$ & $d_2$ & $c_2$ & $b_1$ & $a_1$ & $d_1$ & $c_1$ & $b_3$ & $a_3$ & $d_3$ & $c_3$ \\

$c_2$ & $c_2$ & $d_2$ & $a_2$ & $b_2$ & $c_1$ & $d_1$ & $a_1$ & $b_1$ & $c_3$ & $d_3$ & $a_3$ & $b_3$ \\

$d_2$ & $d_2$ & $c_2$ & $b_2$ & $a_2$ & $d_1$ & $c_1$ & $b_1$ & $a_1$ & $d_3$ & $c_3$ & $b_3$ & $a_3$ \\

$a_3$ & $a_1$ & $b_1$ & $c_1$ & $d_1$ & $a_3$ & $b_3$ & $c_3$ & $d_3$ & $a_2$ & $b_2$ & $c_2$ & $d_2$ \\

$b_3$ & $b_1$ & $a_1$ & $d_1$ & $c_1$ & $b_3$ & $a_3$ & $d_3$ & $c_3$ & $b_2$ & $a_2$ & $d_2$ & $c_2$ \\

$c_3$ & $c_1$ & $d_1$ & $a_1$ & $b_1$ & $c_3$ & $d_3$ & $a_3$ & $b_3$ & $c_2$ & $d_2$ & $a_2$ & $b_2$ \\

$d_3$ & $d_1$ & $c_1$ & $b_1$ & $a_1$ & $d_3$ & $c_3$ & $b_3$ & $a_3$ & $d_2$ & $c_2$ & $b_2$ & $a_2$ \\

\hline

\end{tabular}
\end{center}
{\color{white}.}\\

Denote as $*$ the $1$-block operation induced  by $\bullet$ on
$\GZ$. Let $\Sgg\subset\GZ$ be the topological Markov chain defined
by the oriented graph of Figure \ref{Fig:grafo1}. We have that
$(\Sgg,*)$ is an irreducible quasi-group shift.

\begin{figure}[h]
\centering
\includegraphics[width=0.5\linewidth=1.0]{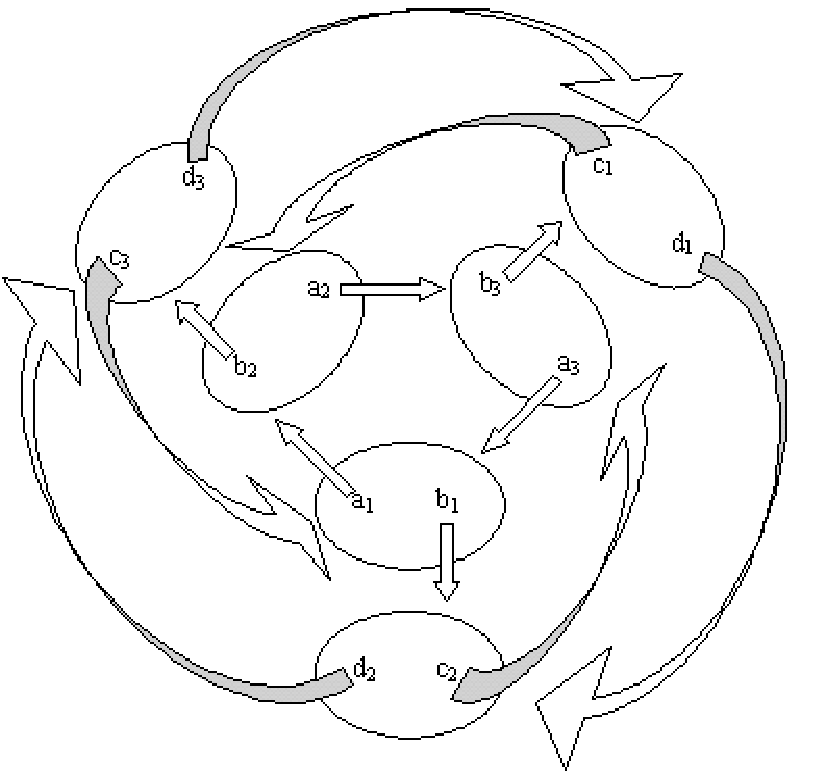}
\caption{Graph generates $\Sgg$. }\label{Fig:grafo1}
\end{figure}

Define the bipermutative cellular automaton $(\Sgg,\Phi)$, where
$\Phi:=id*\s$. It follows $\Phi$ verifies the property (\ref{IV2})
and so it is structurally compatible. Moreover, since $\bullet$ has
the medial property:

$$\forall a,b,c,d\in G,\qquad (a\bullet b)\bullet (c\bullet d)=
(a\bullet c)\bullet (b\bullet d),$$\\

it follows, from (\cite{denes}, Theorem 2.2.2, p.70), that there
exist an Abelian group operation on $G$, $\eta$ and $\rho$ commuting
automorphisms, and $c\in G$, such that $a\bullet
b=\eta(a)+\rho(b)+c$. Therefore, $(\Sgg,\Phi)$ is an affine
c.a..}\\

The next proposition gives a characterization of SC bipermutative
cellular automata.

{\prop\label{IV4}  Let $(\Sgg,\Phi)$ be a SC bipermutative c.a..
Then,

\begin{description}
  \item[(i)] $(\Sgg,\Phi)$ is topologically conjugate to $(\h,\Phi_\h)$ through a
  $1$-block code,
  where $\h=\mathbb{F}\times\Sigma_n$, $\mathbb{F}$ is finite, $\Sigma_n$
  is a {\em full n shift}, and $\Phi_\h=id_\h\otimes\s_\h$ where $\otimes$ is a $k$-block quasi-group
  operation on $\h$.

  \item[(ii)] $\mathbf{h}(\Sgg)=0$ (the topological entropy of the shift is
  zero) if and only if $\Sigma_n=\set{(\ldots,a,a,a,\ldots)}$ (that is, the full shift is
  trivial).

  \item[(iii)] $\Sgg$ is irreducible and has constant sequence if and only if
  $\mathbb{F}=\set{e}$ (that is, $\mathbb{F}$ is unitary).\\
\end{description}}

  \proof {\color{white}.}

  \begin{description}

  \item[(i)] Let $(\Sgg,\Phi)$ be a bipermutative c.a. with radio $1$ which verifies (\ref{IV2}).
  As before, for $a,b\in G$, denote $a\bullet b=\phi(a,b)$ which
  is a quasi-group operation on $G$. Thus,  $\Phi=id*\s$, where
  $*$ is the componentwise quasi-group operation on $\Sgg$ induced from
  $\bullet$.

  From Theorem 4.25 and Remark 4.28 in \cite{Sobottka}, the quasi group $(\Sgg,*)$ is isomorphic
  to a quasi group
  $(\mathbb{F}\times\Sigma_n,\otimes)$, where $\otimes$ is a
  $k$-block operation, with anticipation $k-1$. We denote
  $\h:=\mathbb{F}\times\Sigma_n$, as $\varphi:\Sgg\to \h$ the
  isomorphism between the quasi groups, and $\Phi_\h:=id_{\h}\otimes\s_{\h}$. It follows that

\begin{equation*}
\begin{array}{rll}
\varphi\circ\Phi &=&\varphi\circ
(id_{\Sgg}*\s_{\Sgg})=_{(a)}(\varphi\circ
id_{\Sgg})\otimes(\varphi\circ\s_{\Sgg})\\
&=_{(b)}& (id_\h\circ\varphi)\otimes(\s_{\h}
\circ\varphi)=(id_\h\otimes\s_{\h})
\circ\varphi=\Phi_\h\circ\varphi,
\end{array}
\end{equation*}

where $=_{(a)}$ comes from the fact that $\varphi$ is an isomorphism
between  $(\Sgg,*)$ and $(\h,\otimes)$, and $=_{(b)}$ is due the
fact that $\varphi$ is a 1-block code (see Theorem 4.25 in
\cite{Sobottka}) and so it commutes with the shift map.\\

Since $\otimes$ is a $k$-block quasi-group operation (with memory
$0$ and anticipation $k-1$), we have that $\Phi_\h$ has radio
$k$.\\

\item[(ii) and (iii)] They follow straightforward from Theorem 4.25 of \cite{Sobottka}.

\end{description}

\cqd

{\rem\label{IV5} From Theorem 4.25 in \cite{Sobottka} we could get
an analogous result, but with $\otimes$ being an operation with
memory $k-1$ and anticipation $0$. Therefore, $\Phi_\h$ would
have memory $k$ and anticipation $0$.}\\

We notice $(\h,\Phi_\h)$ in the previous theorem is not necessarily
a bipermutative c.a.. For instance, if $(\Sgg,\Phi)$ is a group c.a.
(which means $(G,\bullet)$ is a group) such that (\ref{IV2}) holds,
then it verifies all hypothesis of Theorem \ref{IV4}, but $\otimes$
can be a $k$-block group operation with memory $0$ and anticipation
$k-1$, for some $k>1$. Thus, $(\h,\Phi_\h)$ can not be right
permutative. In fact, since $\rho:\h_k\times\h_k\to \h_1$ the local
rule of $\otimes$, then since $\s_\h$ is an automorphism to the
group $(\h,\otimes)$, it follows the identity element $\e$ of the
group is such that $\s_\h(\e)=\e$, which implies $\e=(\ldots e,
e,e,\ldots)$. Therefore, taking $w\in\h_k$,
$w=(\overset{k}{\overbrace{e,e,\ldots, e}})$, we have for all $a\in
\h_1$:
$$\phi_\h(wa)=\phi_\h\bigl(
(\overset{k+1}{\overbrace{e,e,\ldots, e,a}})\bigr)=\rho
\bigl((\overset{k}{\overbrace{e,e,\ldots,
e}}),(\overset{k}{\overbrace{e,e,\ldots, a}})\bigr)=e.$$

\section{Right-permutative cellular automata}\label{4}

In this section we shall study two types of cellular automata: {\em
N-scaling}; and {\em $\Psi$-associative}.

We say a cellular automaton $(\Sgg,\Phi)$ with radio $1$ is a
N-scaling c.a. for some $N\geq 2$ if its local rule $\phi:G\times
G\to G$ is such that for any $\mathbf{x}=(x_i)_{i\in\Z}\in\Sgg$,
$$\bigl(\Phi^N(\mathbf{x})\bigr)_0=x_0\bullet x_N.$$

On the other hand $(\Sgg,\Phi)$ is said $\Psi$-associative, if there
exists a permutation $\Psi:G\to G$ such that for any $a,b,c\in G$,
we have
$$(a\bullet b)\bullet c = \Psi\bigl(a\bullet(b\bullet c)\bigr)$$

When $\Sgg=\GZ$, Host-Maass-Mart\'{\i}nez \cite{HostMaassMartinez} have
proved that every right-permutative\sloppy\ N-scaling c.a.
$(\Sgg,\Phi)$ is topologically conjugate to the product of an affine
c.a. with a translation, while every right-permutative
$\Psi$-associative is topologically conjugate to the product of a
group c.a. with a translation.

Theorems \ref{IV7a} and \ref{IV7b} below reproduce those results for
the general case of cellular automata defined on topological Markov
chains. To proof these theorems we shall remark some basics on these
types of cellular automata:

    {\rem\label{IV6}  {\color{white} . }
\begin{itemize}

    \item If $(\GZ,\Phi)$ is a right-permutative $\psi$-associative c.a., from Theorem 6 in
    \cite{HostMaassMartinez},
    we get that there exists a $1$-block code $\mathbf{u}:\GZ\to K^\Z\times B^\Z$\sloppy,
    which is a topological conjugacy between
    $(\GZ,\Phi)$ and $(K^\Z\times B^\Z,\Phi_K\times \mathbf{g}_B)$,
    where $B\subseteq G$ and $K$ are two finite alphabets, $\phi_K$ is
    a {\em group c.a.} and $\mathbf{g}_B$ is a
    translation.

    We recall $\mathbf{g}_B=\mathbf{s}_B\circ\s_B$,
    where $\mathbf{s}_B:B^\Z\to B^\Z$ is a 1-block code with
    local rule $s_B:B\to B$ which is a permutation on $B$.
    Moreover, \cite{HostMaassMartinez} gives $s_B:B\to B$ is defined
    for all $e'\in B$ by $s_B(e')=e''\bullet e'$, where $e''\in B$ is
    any element.

    Furthermore, $\mathbf{u}$ has local rule $u:G\to K\times
    B$ which is a bijection and is given for any $a\in G$ by
    $$u(a)=(\tilde{a},e_a),$$

    where $\tilde{a}$ is the equivalent class of $a\in G$ to the equivalence relation,
    $$a\sim b\Longleftrightarrow \forall c\in G,\ a\bullet c=b\bullet
    c,$$

    and $e_a$ is the unique
    element of $G$ for which
    $a\bullet e_a=a$.
    We notice for all $a\in G$ and $e\in B$, we have
    $a\bullet e\sim a$. Moreover, the following property holds:
    $e_{a\bullet b}=e_a\bullet e_b=s_B(e_b)$.

    Finally, since $\Phi_K$ is a group, its local rule
    define a group operation on $K$:
    $$\forall\tilde{a},\tilde{b}\in K,
    \tilde{a}\tilde{\bullet}\tilde{b}:=\phi_K(\tilde{a},\tilde{b}).$$\\

    \item From Theorem 8 in \cite{HostMaassMartinez},
    if $(\GZ,\Phi)$ is a $N$-scaling c.a, then the above statements hold, but $\Phi_K$
    will be an affine c.a. and in the code $u(a)=(\tilde{a},e_a)$, $e_a$
    will be defined as the unique element of $B$ for which the equation $e_a=x\bullet a$
    has solution.\\

\end{itemize}}

{\theo\label{IV7a}  Let $(\Sgg,\Phi)$ be a SC right-permutative
$\Psi$-associative c.a.. Then, $(\Sgg,\Phi)$ is topologically
conjugate through a $1$-block code to
$(\mathfrak{K}\times\mathfrak{B},\Phi_\mathfrak{K}\times\mathbf{g}_\mathfrak{B})$,
where $\mathfrak{K}$ and $\mathfrak{B}$ are topological Markov
chains, $(\mathfrak{K},\Phi_\mathfrak{K})$ is a group c.a., and
$(\mathfrak{B},\mathbf{g}_\mathfrak{B})$ is a translation.}

 \proof{\color{white}.}

\begin{description}

  \item[Step 1] Since $(\Sgg,\Phi)$ has radio $1$ and verifies (\ref{IV2})
  we can consider that $\Phi:\Sgg\to\Sgg$ is a
  restriction on $\Sgg$ of some right-permutative $\Psi$-associative c.a.
  $(\GZ,\tilde{\Phi})$ which has the same local rule $\phi:G\times G\to G$.

  Let $(K^\Z\times B^\Z,\Phi_K\times \mathbf{g}_B)$ be the
  group-translation and $\mathbf{u}:\GZ\to K^\Z\times B^\Z$ the
  topological conjugacy presented in Remark \ref{IV6}.

  We consider on $K\times B$ the
  right-permutative operation also denoted as $\bullet$ and
  induced from the local rule of
  $\Phi_K\times \mathbf{g}_B$: given
  $(\tilde{a}_1,e_1),(\tilde{a}_2,e_2)\in K\times B$ define

  $$(\tilde{a}_1,e_1)\bullet (\tilde{a}_2,e_2)=\bigl(\tilde{a}_1\tilde{\bullet}
  \tilde{a}_2,s_B(e_2)\bigr).$$

  Notice that $u:G\to K\times B$ is an isomorphism between
  $(G,\bullet)$ and $(K\times B,\bullet)$. In fact, $u$ is
  bijective and

  $$u(a\bullet c)=(\widetilde{a\bullet c},e_{a\bullet c})=_{(a)}
  (\tilde{a}\tilde{\bullet}\tilde{c},e_a\bullet e_c)$$
  $$=\bigl(\tilde{a}\tilde{\bullet}\tilde{c},s_B(e_c)\bigr)=
  (\tilde{a},e_a)\bullet (\tilde{c},e_c)=u(a)\bullet u(c),$$

  where $=_{(a)}$ comes from Theorem 6 of \cite{HostMaassMartinez}.

  The operation $\bullet$ on $K\times B$ induces
  the componentwise operation also denoted as $*$ on $K^\Z\times B^\Z$. Thus,
  $\mathbf{u}:\GZ\to K^\Z\times B^\Z$ is an isomorphism between
  $(\GZ,*)$ and $( K^\Z\times B^\Z,*)$.

  Define $\Lambda:=\mathbf{u}(\Sgg)\subseteq K^\Z\times B^\Z$.
  Since $\mathbf{u}$ is topological conjugacy between $\Phi$
  and $\Phi_K\times\mathbf{g}_B$, it follows
  $\Phi_K\times\mathbf{g}_B(\Lambda)=\Lambda$. Therefore, we have
  the cellular automaton $(\Lambda,\Phi_K\times\mathbf{g}_B)$ is
  well defined and $*$ is closed on $\Lambda$. Moreover, $\mathbf{u}|_{\Sgg}$ is
  a topological conjugacy between $(\Sgg,\Phi)$ and
  $(\Lambda,\Phi_K\times\mathbf{g}_B)$, and an isomorphism between $(\Sgg,*)$
  and $(\Lambda,*)$.\\

  \item[Step 2] We will show that there exists $M\geq 1$ such that for all
  $e\in B$ we have $s_B^{M}(e)=e$.

  Since $s_B$ is a permutation on $B$, it follows for all $e\in B$
  there exists $M_{e}\geq 1$ such that $s_B^{M_{e}}(e)=e$. Because
  $B$ is a finite alphabet, we can take $M$  a multiple of all
  periods of each element of $B$. Then, the result follows.\\

  \item[Step 3]  Let us to prove that
  $\Lambda=\mathfrak{K}\times\mathfrak{B}$, where
  $\mathfrak{K}\subseteq K^\Z$ and $\mathfrak{B}\subseteq B^\Z$ are
  both topological Markov chains.

  First, notice that, because $\tilde{\bullet}$ is a quasi-group
  operation,
  there exists $L\in\N$ such that for all $\tilde{a},\tilde{c}\in
  K$:

  $$\underbrace{(\ (\ (\tilde{c}\tilde{\bullet}
  \tilde{a})\ldots\tilde{\bullet} \tilde{a})\tilde{\bullet}
  \tilde{a})\tilde{\bullet} \tilde{a}}_{\tilde{c}\ multiplied\ L\
  times\ by\ \tilde{a}\ for\ the\ right\ side}
  =\tilde{c}.$$\\

  Denote $\pi_K:\Lambda\to K^\Z$ and $\pi_B:\Lambda\to B^\Z$ the
  canonical projections on the first and second coordinates
  respectively.

  It is straightforward that $\Lambda\subseteq
  \pi_K(\Lambda)\times\pi_B(\Lambda)$. So, we only need to show that
  $\pi_K(\Lambda)\times\pi_B(\Lambda)\subseteq \Lambda$.

  In fact, given $(\tilde{c}_i)_{i\in\Z}\in\pi_K(\Lambda)$ and
  $(e_i)_{i\in\Z}\in\pi_B(\Lambda)$ must there exist
  $(\tilde{c}_i,e'_i)_{i\in\Z},(\tilde{a}_i,e_i)_{i\in\Z}\in\Lambda$,
  and so

  $$ \underbrace{\Bigl(\Bigl(\Bigl((\tilde{c}_i,e'_i)_{i\in\Z}*(\tilde{a}_i,e_i)_{i\in\Z}\Bigr)*
  (\tilde{a}_i,e_i)_{i\in\Z}\Bigr)\ldots*(\tilde{a}_i,e_i)_{i\in\Z}\Bigr)
  *(\tilde{a}_i,e_i)_{i\in\Z}}_{multiplying\
  (\tilde{c}_i,e'_i)_{i\in\Z}\ L\ times\ by\
  (\tilde{a}_i,e_i)_{i\in\Z}\ for\ the\ right\ side}$$

  $$=\Bigl(\underbrace{(((\tilde{c}_i\tilde{\bullet}\tilde{a}_i)\tilde{\bullet}
  \tilde{a}_i)\ldots\tilde{\bullet}\tilde{a}_i)\tilde{\bullet}
  \tilde{a}_i}_{multiplying\ \tilde{c}_i\ L\ times\ by\ \tilde{a}_i\
  for\ the\ right\ side},s_B(e_i)\Bigr)_{i\in\Z}$$

  $$=\bigl(\tilde{c}_i,s_B(e_i)\bigr)_{i\in\Z}\in\Lambda.$$\\

  Now, we repeat the above procedure, but multiplying
  $\bigl(\tilde{c}_i,s_B(e_i)\bigr)_{i\in\Z}$ $L$ times by itself
  for the right side, and so we get
  $\bigl(\tilde{c}_i,s_B^2(e_i)\bigr)_{i\in\Z}\in\Lambda$. By
  induction, we can obtain that for all $m\geq 0$,
  $\bigl(\tilde{c}_i,s_B^m(e_i)\bigr)_{i\in\Z}\in\Lambda$. From Step
  2 there exists $M\geq 1$ such that for all $i\in\Z$ we have
  $s_B^{M}(e_i)=e_i$. Therefore, we get
  $\bigl(\tilde{c}_i,e_i\bigr)_{i\in\Z}\in\Lambda$, which allows us
  to deduce that $\Lambda=\pi_K(\Lambda)\times\pi_B(\Lambda)$.

  Notice that $\mathbf{u}|_{\Sgg}$ is a $1$-block code from $\Sgg$
  to $\Lambda$ such that its inverse is also a $1$-block code. Thus,
  since $\Sgg$ is a topological Markov chain, it follows that
  $\Lambda$ is also a topological Markov chain.
  Finally, since
  $\Lambda=\pi_K(\Lambda)\times\pi_B(\Lambda)$ we have that
  $\pi_K(\Lambda)$ and $\pi_B(\Lambda)$ are also both topological Markov
  chains, and denoting $\mathfrak{K}:=\pi_K(\Lambda)$,
  $\mathfrak{B}:=\pi_B(\Lambda)$,
  $\Phi_\mathfrak{K}:=\Phi_K|_\mathfrak{K}$ and
  $\Phi_\mathfrak{B}:=\Phi_B|_\mathfrak{B}$ we finish the proof.

\end{description}

\cqd

{\theo\label{IV7b}  Let $(\Sgg,\Phi)$ be a SC right-permutative
$N$-scaling c.a.. If its extension $(\GZ,\tilde{\Phi})$ is also a
$N$-scaling c.a., then $(\Sgg,\Phi)$ is topologically conjugate
through a $1$-block code to
$(\mathfrak{K}\times\mathfrak{B},\Phi_\mathfrak{K}\times\mathbf{g}_\mathfrak{B})$,
where $\mathfrak{K}$ and $\mathfrak{B}$ are topological Markov
chains, $(\mathfrak{K},\Phi_\mathfrak{K})$ is an affine c.a., and
$(\mathfrak{B},\mathbf{g}_\mathfrak{B})$ is a translation.}

\proof

Since $(\Sgg,\Phi)$ is the restriction on $\Sgg$ of a $N$-scaling
c.a. $(\GZ,\tilde{\Phi})$, we can apply a similar reasoning than
Theorem \ref{IV7a}.

\cqd

{\cor\label{IV7c}  Let $(\Sgg,\Phi)$ be a SC right-permutative
$N$-scaling c.a.. If $\Sgg$ is mixing, then $(\Sgg,\Phi)$ is
topologically conjugate through a $1$-block code to
$(\mathfrak{K}\times\mathfrak{B},\Phi_\mathfrak{K}\times\mathbf{g}_\mathfrak{B})$,
where $\mathfrak{K}$ and $\mathfrak{B}$ are topological Markov
chains, $(\mathfrak{K},\Phi_\mathfrak{K})$ is an affine c.a., and
$(\mathfrak{B},\mathbf{g}_\mathfrak{B})$ is a translation.}

\proof

Since $\Sgg$ is mixing, there exists $q\geq 1$ such that for any
$k\geq q$ and $u,w\in G$ we always can find
$(v_1,\ldots,v_k)\in\Sgg_k$ such that
$(u,v_1,\ldots,v_k,w)\in\Sgg_{k+2}$. Without loss of generality we
can consider $N\geq q$, because if $(\Sgg,\Phi)$ is $N$-scaling,
then it is also $N^m$-scaling for any $m\geq 1$.
We will show that $(\GZ,\tilde{\Phi})$ is also $N$-scaling:\\

Given a sequence $\mathbf{x}=(x_i)_{i\in\Z}\in G$, due the fact of
$\Sgg$ is mixing and $N\geq q$, we can find a sequence
$\mathbf{y}=(y_i)_{i\in\Z}\in\Sgg$ such that $y_{jN}=x_{j}$ for all
$i\in\Z$. Thus,

$$\bigl(\tilde{\Phi}(\mathbf{x})\bigr)_{j}=x_j\bullet x_{j+1}=y_{jN}\bullet y_{(j+1)N}
=\bigl(\Phi^N(\mathbf{y})\bigr)_{jN}\ ,$$

and by induction we get that for any $k\geq 1$,
$\bigl(\tilde{\Phi}^k(\mathbf{x})\bigr)_{j}=\bigl(\Phi^{kN}(\mathbf{y})\bigr)_{jN}$.
Therefore,
$$\bigl(\tilde{\Phi}^{N}(\mathbf{x})\bigr)_{j}=\bigl(\Phi^{N^2}(\mathbf{y})\bigr)_{jN}=y_{jN}\bullet y_{jN+N^2}
=y_{jN}\bullet y_{(j+N)N}=x_j\bullet x_{j+N} \ .$$

Now, since $\tilde{\Phi}$ is a $N$-scaling c.a., we can apply
Theorem \ref{IV7b} to conclude the proof.

\cqd

Notice that $(\mathfrak{K},\Phi_K)$ obtained in the previous
theorems is a group c.a. (or an affine c.a.) which is also
structurally compatible. Thus, since $(\mathfrak{K},\Phi_K)$ is
bipermutative, we can apply Proposition \ref{IV4} to get it is
topologically conjugate through a $1$-block code to $(\h,\Phi_\h)$,
where $G=\mathbb{F}\times\Sigma_n$ with $\mathbb{F}$ is finite and
$\Sigma_n$ is a full n shift, and $\Phi_G$ is a group c.a. (or an
affine c.a.) with radio $k$.


\section{Projections of measures with complete connections and summable
decay}\label{5}

In this section we shall present sufficient conditions to reproduce
results about the convergence of the Ces\`{a}ro mean distribution
(\cite{HostMaassMartinez}, \cite{mms}) to the more general case of
$\Sgg$ being neither a full shift nor a groupshift, but
$(\Sgg,\Phi)$ being structurally compatible.

{\lem Let $\Lambda$ and $\Lambda'$ be two topological Markov chains,
and $\Theta:\Lambda\to\Lambda'$ be an invertible $1$-block code
which action is constant on the predecessor sets. Suppose
$\Theta^{-1}$ has memory $1$ and anticipation $0$. If $\mu$ is a
$\s$-invariant probability measure on $\Lambda$ with complete
connections (compatible with $\Lambda$) and summable decay, then
$\mu'=\mu\circ\Theta^{-1}$ also has complete connections (compatible
with $\Lambda'$) and summable decay.}

\proof

Let $C'$ be a cylinder of $\Lambda'$ defined by the coordinates
$i=0,\ldots,m$ with $m\geq 1$, that is, $C'=[c_0',\ldots,c_m']$. We
will show that $C:=\Theta^{-1}(C')$ is a cylinder of $\Lambda$
defined by the coordinates $i=1,\ldots,m$, that is,
$C=[c_1,\ldots,c_m]$.

Denote as $\theta$ the local rule of $\Theta$. Notice that for all
$1\leq i\leq m$, $c_i\in \Lambda_1$ is well defined by
$c_i:=\theta\n(c_{i-1}',c_i')$. Therefore,

$$\Theta\n(C')=\bigcup_{\left.
\begin{array}{l} c_0\in\Pp(c_1)\\
\theta(c_0)=c_0'\end{array}\right.}[c_0,c_1,\ldots,c_m]
=\bigcup_{c_0\in\Pp(c_1)}[c_0,c_1,\ldots,c_m]=[c_1,\ldots,c_m]$$\\

Through the use of a similar reasoning and since $\Theta\n$ has
anticipation $0$, we get that for any $v',w'\in\Lambda'^-$, we can
define $v:=\Theta\n(v')$ and $w:=\Theta\n(w')$ which are both pasts
belonging to $\Lambda^-$. In particular, if $v_{-i}'=w_{-i}'$ for
$1\leq i\leq m$, with $m\geq 2$, then $v_{-i}=w_{-i}$ for $1\leq
i\leq m-1$.

On the other hand, since $\mu$ has complete connections (compatible
with $\Lambda$), given $w'\in\Lambda'^-$ and $a'\in\F(w_{-1}')$
there exist unique $w\in\Lambda$ and $a\in\F(w_{-1})$ such that
$\mu'_{w'}(a')=\mu_w(a)>0$. It means $\mu'$ also has complete
connections (compatible with $\Lambda'$). Moreover, for $m\geq 2$,
it follows

$$\gamma_m'=\sup\set{\abs{\frac{\mu'_{v'}(a')}{\mu'_{w'}(a')}-1}:\ \begin{array}{l}v',w'\in
\Lambda'^-;\quad v'_{-i}=w'_{-i},\ 1\leq i\leq m;\\
a'\in\F(v'_{-1})=\F(w'_{-1})\end{array}}$$\\

$$=\sup\set{\abs{\frac{\mu_v(a)}{\mu_w(a)}-1}:\ \begin{array}{l}v,w\in
\Lambda^-;\quad v_{-i}=w_{-i},\ 1\leq i\leq m-1;\\
a\in\F(v_{-1})=\F(w_{-1})\end{array}}=\gamma_{m-1},$$\\

which means $\mu'$ has summable decay.

\cqd

In an analogous way, we can prove the following Lemma:

{\lem\label{IV10} Let $\Lambda$ and $\Lambda'$ be two topological
Markov chains, and let
$\varphi:\Lambda\times\Sigma\to\Lambda'\times\Sigma$ be block code
defined by $\varphi:=\Theta\times id$, where
$\Theta:\Lambda\to\Lambda'$ is an invertible $1$-block code which is
constant on the predecessor sets. Suppose $\Theta^{-1}$ has memory
$1$ and anticipation $0$. If $\mu$ is $\s$-invariant probability
measure on $\Lambda$ with complete connections (compatible with
$\Lambda$) and summable decay, then $\mu'=\mu\circ\Theta^{-1}$ also
has complete connections (compatible with $\Lambda'$) and summable
decay.}

\cqd

Now, consider $(\Sgg,\Phi)$ being a SC bipermutative c.a.. Let
$\varphi:\Sgg\to G$ be the topological conjugacy between
$(\Sgg,\Phi)$ and $(\h,\Phi_\h)$, where
$\h=\mathbb{F}\times\Sigma_n$, given by Proposition \ref{IV4}. From
Remark \ref{IV5} we can suppose that $\varphi$ has memoria $k$ and
anticipation $0$. With this notations, we have that:

{\prop\label{IV11} If $(\Sgg,\Phi)$ is a SC bipermutative c.a., and
$\mu$ is a probability measure with complete connections (compatible
with $\Sgg$) and summable decay, then $\mu\circ\varphi^{-1}$ is a
probability measure on $\h=\mathbb{F}\times\Sigma_n$ which also has
complete connections and summable decay.}

\proof

From Theorem 4.25 of \cite{Sobottka}, $\varphi$ is given by the
following composition:
$$\varphi=\varphi_n\circ\eta_n\circ\varphi_{n-1}\circ\eta_{n-1}\circ\ldots\circ\varphi_1\circ\eta_1,$$
where for all $i=1,\ldots, n$, $\varphi_i=\Theta_i\times id$ is a
block code as in Lemma \ref{IV10}, and $\eta_i$ is an invertible
$1$-block code which inverse is also a $1$-block code. Thus, for
each $i\leq n$ we have that $\eta_i$ and $\varphi_i$ preserve the
properties of complete connections and summable decay of the
measure, which conclude the proof.

\cqd

\section{Ces\`{a}ro mean convergence of measures with complete connections and summable
decay}\label{6}

In this section we shall present some results about the convergence
of the Ces\`{a}ro mean distribution of probability measures under the
action cellular automata, namely we study the following limit:
$$\lim_{N\rightarrow\infty}\frac{1}{N}\sum_{n=0}^{N-1}\mu\circ\Phi^{-n}.$$

The essential tools that we will use to study the convergence of the
Ces\`{a}ro mean distribution are: propositions \ref{IV4} and \ref{IV11};
and Corollary 29 of \cite{HostMaassMartinez}.

{\defn Given a SC bipermutative c.a. $(\Sgg,\Phi)$, we say it is
{\em regular} if the quasigroup $(\mathbb{F}\times\Sigma_n,\otimes)$
set in Theorem \ref{IV4} is such that
$\otimes=\otimes_\mathbb{F}\times\otimes_{\Sigma_n}$, where
$(\mathbb{F},\otimes_\mathbb{F})$ and
$(\Sigma_n,\otimes_{\Sigma_n})$ are both quasigroups. Furthermore,
if $\otimes_{\Sigma_n}$ is a $1$-block operation, then we say
$(\Sgg,\Phi)$ is {\em simple}.}\\

{\ex If $\Sgg$ is irreducible or $\mathbf{h}(\Sgg) = 0$, then $\Sgg$
is regular due Theorem 4.25(ii,iii) of \cite{Sobottka}. If
$\mathbf{h}(\Sgg) = p$, where $p$ is a prime number, then $\Sgg$ is
simple due Theorem 4.26 of \cite{Sobottka}.}\\

{\theo\label{IV12} Let $(\Sgg,\Phi)$ be a SC cellular automaton,
where $\Sgg$ is not necessarily irreducible. Denote as
$(\GZ,\tilde{\Phi})$ the extension of $(\Sgg,\Phi)$ to the full
shift, and suppose $\mu$ is a probability measure on $\Sgg$ with
complete connections (compatible with $\Sgg$) and summable decay.
Then:
\begin{description}

\item[(i)] If $(\Sgg,\Phi)$ is an affine c.a. which is regular and simple,
then the Ces\`{a}ro mean distribution of $\mu$ under the action of
$\Phi$ converges to a maximum entropy measure. In particular, if
$\Sgg$ is irreducible and has a constant sequence, then the Ces\`{a}ro
mean distribution of $\mu$ under the action of $\Phi$ converges to
the unique maximum entropy measure (the Parry measure);

\item[(ii)] If $(\Sgg,\Phi)$ is a right-permutative
$\Psi$-associative c.a. and the group c.a. associate to it (see
Theorem \ref{IV7a}) is Abelian, regular and simple, then the Ces\`{a}ro
mean distribution of $\mu$ under the action of $\Phi$ converges.

\item[(iii)] If $(\GZ,\tilde{\Phi})$ is right-permutative and $N$-scaling
and the affine c.a. associate to it (see Theorem \ref{IV7b}) is
regular and simple, then the Ces\`{a}ro mean distribution of $\mu$ under
the action of $\Phi$ converges.

\end{description}}

\proof{\color{white}.}

\begin{description}

\item[(i)] Let
$(\mathbb{F}\times\Sigma_n,\Phi_{\mathbb{F}\times\Sigma_n})$ and
$\varphi:\Sgg\to\mathbb{F}\times\Sigma_n$ be the cellular automaton
and the topological conjugacy given by Proposition \ref{IV4}. From
Proposition \ref{IV11}, we have that $\mu'=\mu\circ\varphi^{-1}$ is
a probability measure on $\mathbb{F}\times\Sigma_n$ with complete
connections and summable decay. Moreover, since $(\Sgg,\Phi)$ is
regular and simple, it follows that
$\Phi_{\mathbb{F}\times\Sigma_n}=\Phi_\mathbb{F}\times\Phi_{\Sigma_n}$,
where 
$(\Sigma_n,\Phi_{\Sigma_n})$ is an affine c.a.. In fact,
$(\mathbb{F}\times\Sigma_n,\otimes_{\mathbb{F}\times\Sigma_n})=
(\mathbb{F}\times\Sigma_n,\otimes_\mathbb{F}\times\otimes_{\Sigma_n})$
has the medial property, thus 
$(\Sigma_n,\otimes_{\Sigma_n})$ also has the medial property and we
can apply (\cite{denes}, Theorem 2.2.2, p.70), in the same way as in
Theorem \ref{IV8}, which allows us to deduce that $\Phi_{\Sigma_n}$
is an affine c.a..

Furthermore, since $\mathbb{F}$ is a finite set, we get
$(\mathbb{F},\Phi_{\mathbb{F}})$ is equicontinuous. Therefore, from
Corollary 29 in \cite{HostMaassMartinez}, it follows that the Ces\`{a}ro
mean of $\mu'$ under the action of $\Phi_{\mathbb{F}\times\Sigma_n}$
converges to a probability measure $\mu'_{\mathbb{F}}\times\nu$,
where $\mu'_{\mathbb{F}}$ is a $\Phi_{\mathbb{F}}$-invariant
probability measure on $\mathbb{F}$, and $\nu$ is the Parry measure
on $\Sigma_n$ (that is, the uniform Bernoulli measure). Since
$(\Sgg,\Phi)$ is topologically conjugate to
$(\mathbb{F}\times\Sigma_n,\Phi_{\mathbb{F}\times\Sigma_n})$, we
conclude that
\begin{equation}\label{limit}\lim_{N\rightarrow\infty}\frac{1}{N}\sum_{n=0}^{N-1}\mu\circ\Phi^{-n}=(
\mu'_{\mathbb{F}}\times\nu)\circ\varphi,\end{equation}\\

which is a maximum entropy measure since it is the projection of a
measure on a finite set $\mathbb{F}$ product the Parry measure on
the full n-shift.

In particular, when $\mathbb{F}$ is not unitary (and since it is
finite and hence has zero entropy) we conclude there could exist
more than one maximum entropy measure for $\Phi$ (\eqref{limit} is
one of them). On the other hand, if $\Sgg$ is irreducible and has a
constant sequence, then, from Proposition \ref{IV4}(iii),
$\mathbb{F}$ is unitary and the limit measure is exactly the
projection of the Parry measure on the full n-shift. In such
case the limit measure is the Parry measure on $\Sgg$.\\

\item[(ii)]

From Theorem \ref{IV7a} and Proposition \ref{IV4}, and using
(\cite{denes}, Theorem 2.2.2, p.70) in the same way as in the proof
of Theorem \ref{IV8}, we deduce that $(\Sgg,\Phi)$  can be
represented as $(\mathfrak{B}\times\mathbb{F}\times
\Sigma_n,\mathbf{g}_B\times\Phi_{\mathbb{F}}\times\Phi_{\Sigma_n})$,
that is: a translation on a topological Markov chain, product a
group c.a. on a finite set, product a group c.a. on a full shift.
Thus, by Corollary 29 of \cite{HostMaassMartinez}, we conclude
\begin{equation}\lim_{N\rightarrow\infty}\frac{1}{N}\sum_{n=0}^{N-1}\mu\circ\Phi^{-n}=
\mu_{\mathfrak{B}}\times(\mu'_{\mathbb{F}}\times\nu)\circ\varphi,\end{equation}\\

where $\mu_{\mathfrak{B}}$ is a $\mathbf{g}_B$-invariant probability
measure on $\mathfrak{B}$, $\mu'_{\mathbb{F}}$ is a
$\Phi_{\mathbb{F}}$-invariant probability measure on $\mathbb{F}$,
and $\nu$ is the Parry measure on $\Sigma_n$ (that is, the uniform
Bernoulli measure).\\

\item[(iii)] This proof is analogous to the part (ii), but uses Theorem
\ref{IV7b} instead Theorem \ref{IV7a}.

\end{description}

\cqd

{\ex The affine c.a. of Example \ref{IV3} is regular and simple,
because applying the reasoning presented in the proof of Proposition
\ref{IV4} we deduce it is topologically conjugate to
$(\mathbb{F}\times\Z_2^\Z,\Phi_{\mathbb{F}}\times\Phi_{\Z_2^\Z})$,
where $\mathbb{F}=\set{(\ldots,0,1,2,0,1,2,\ldots)}$,
$\Phi_{\mathbb{F}}=\s_\mathbb{F}$ and $\Phi_{\Z_2^\Z}= id+\s$.
Therefore, the Ces\`{a}ro mean of any probability measure on $\Sgg$ with
complete connections summable decay converges under the action of
$(\Sgg,\Phi)$.}\sloppy

{\ex\label{IV13} Let $\bullet$ be a binary operation defined on
the set $G=\set{a,b,c,d,e,f,g,h}$ by the following table:\\

\begin{center}
\begin{tabular}{|c||c c c c c c c c|}
\hline

$\bullet$

    & $a$ & $b$ & $c$ & $d$ & $e$ & $f$ & $g$ & $h$ \\
\hline\hline

$a$ & $b$ & $a$ & $d$ & $c$ & $f$ & $e$ & $h$ & $g$ \\

$b$ & $b$ & $a$ & $d$ & $c$ & $f$ & $e$ & $h$ & $g$ \\

$c$ & $d$ & $c$ & $b$ & $a$ & $h$ & $g$ & $f$ & $e$ \\

$d$ & $d$ & $c$ & $b$ & $a$ & $h$ & $g$ & $f$ & $e$ \\

$e$ & $f$ & $e$ & $h$ & $g$ & $b$ & $a$ & $d$ & $c$ \\

$f$ & $f$ & $e$ & $h$ & $g$ & $b$ & $a$ & $d$ & $c$ \\

$g$ & $h$ & $g$ & $f$ & $e$ & $d$ & $c$ & $b$ & $a$ \\

$h$ & $h$ & $g$ & $f$ & $e$ & $d$ & $c$ & $b$ & $a$ \\

\hline

\end{tabular}\\
\end{center}

Let $\Lambda\subset\GZ$ be the topological Markov chain defined by
the oriented graph presented in Figure \ref{Fig:grafo2}.

\begin{figure}[htbp]
\centering
\includegraphics[width=0.14\linewidth=1.0]{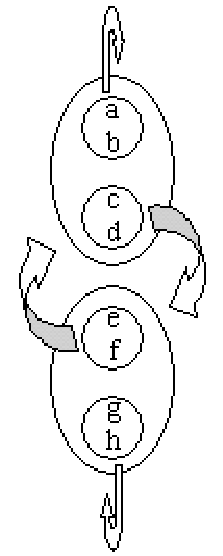}
\caption{Graph which generates $\Lambda$. }\label{Fig:grafo2}
\end{figure}

We define the map $\phi:G\times G\to G$ giving by
$\phi(u,v)=u\bullet v$, and we consider the cellular automaton
$\Phi:\Lambda\to\Lambda$ with radio $1$, which the local rule is
$\phi$.

It is easy to check that $(\Lambda,\Phi)$ is a SC right-permutative
and $\Psi$-associative c.a., where $\Psi$ is given by
$\Psi(\cdot):=\phi(a,\cdot)$.

From the algorithm developed in the proof of Theorem \ref{IV7a}, we
get $(\Lambda,\Phi)$ is topologically conjugate to
$\bigl(\Sigma\times\set{0,1}^\Z,\Phi_\Sigma\times \mathbf{g}\bigr)$,
where: $\Sigma\subset(\Z_2\oplus\Z_2)^\Z$; $\Phi_\Sigma$ is a group
c.a.; and $\mathbf{g}$ is a translation on $\set{0,1}^\Z$ giving by
$\mathbf{g}\bigl((x_i)_{i\in\Z}\bigr)=\bigl(\Psi(x_{i+1})\bigr)_{i\in\Z}$.

Using the algorithm developed in the proof of Proposition \ref{IV4}
we find $(\Sigma,\Phi_\Sigma)$ is topologically conjugate to the
group c.a. $(\Z_2^\Z,\Phi_{\Z_2^\Z})$, and so it is regular and
simple. Therefore, Theorem \ref{IV12} guarantees the convergence of
the Ces\`{a}ro mean of any probability measure on $\Lambda$ with
complete connections summable decay under the action of
$(\Lambda,\Phi)$.}


\section{Invariant measures for cellular automata on topological Markov
chains}\label{7}

Let $(\Sgg,\Phi)$ be a SC cellular automaton and suppose that $\Sgg$
is irreducible.

An important problem is to characterize probability measures on
$\Sgg$ which are invariant for the $\Z^2$-action defined on $\Sgg$
by $(\Phi,\s)$. Several works (\cite{HostMaassMartinez},
\cite{Pivato}, \cite{Sablik}) have  studied this problem and for
many cases have showed that the Parry measure (the unique maximum
entropy measure for $(\Sgg,\s)$) is the unique $(\s,\Phi)$-invariant
measure.

We can deduce results about $(\Phi,\s)$-invariant measures through
the use of the topological conjugacies presented previously. When
$\Sgg$ is a group shift, then the following theorems are particular
cases of the results presented by Sablik \cite{Sablik}.

{\theo\label{IV8} Let $(\Sgg,\Phi)$ be a SC affine c.a, with $\Sgg$
being irreducible and $\mathbf{h}(\Sgg)=\log p$, where $p$ is a
prime number. Let $\mu$ be a $(\Phi,\s)$-invariant probability
measure on $\Sgg$. If $\mu$ is ergodic to $\s$ and has
positive entropy to $\Phi$, then $\mu$ is the Parry measure.}\\

\proof

Since $(\Sgg,\Phi)$ is an affine c.a. it follows that it is
bipermutative, which implies the operation $\bullet$, defined by
$a\bullet b:=\phi(a,b)$ is a quasi-group operation on $G$.

On the other hand, from definition of affine c.a. there exists an
Abelian group operation on $G$, $\eta$ and $\rho$ commuting
automorphism and $k\in G$, such that $\phi(a,b)=\eta(a)+\rho(b)+k$.
It implies that $\bullet$ has the medial property. Thus, the
componentwise quasi-group operation $*$ induced from $\bullet$ on
$\Sgg$, also has the medial property.

From Proposition \ref{IV4}, $(\Sgg,\Phi)$ is topologically conjugate
to $(K^\Z,\Phi_K)$ through a 1-block code, where $\Phi_K$ is given
by $\Phi_K=id\otimes\s$. Moreover the same code is an isomorphism
between $(\Sgg,*)$ and $(K^\Z,\otimes)$. Therefore, $\otimes$ is
also a quasi-group operation which has the
medial property.\\

Since $\mathbf{h}(\Sgg)=\log p$, with $p$ being a prime number, from
Theorem 4.26 of \cite{Sobottka} gives $\abs{K}=p$ and
$\otimes$ is a 1-block operation. Thus, 
there exists a quasi-group operation $\odot$ on $K$, which induces
the operation $\otimes$. Notice that the local rule of $\Phi_K$ is
given by $\phi_K(a',b')=a'\odot b'$.

Hence, $\odot$ also has the medial property, and so from
(\cite{denes}, Theorem 2.2.2, p.70) there exist an Abelian group
operation $\oplus$ on $K$, two commuting automorphism $\eta'$ and
$\rho'$, and $c'\in K$, such that $a'\odot
b'=\eta'(a')\oplus\rho'(b')\oplus c'$. With other words,
$(K^\Z,\Phi_K)$ is an affine c.a..\\

Now, defining $\mu':=\mu\circ\varphi^{-1}$, we have that
$(K^\Z,\mathbf{g})$ and $\mu'$ verify all hypothesis of Theorem 12
in \cite{HostMaassMartinez} which implies $\mu'$ is the uniform
Bernoulli measure on $K^\Z$, i.e., the maximum entropy measure for
the full shift. Therefore, we conclude that $\mu$ is the maximum
entropy measure on $\Sgg$.

\cqd

The following theorem has analogous proof than the previous one, but
uses Theorem 13 instead Theorem 12 of \cite{HostMaassMartinez}.

{\theo Let $(\Sgg,\Phi)$ be a SC affine c.a., such that $\Sgg$ is
irreducible and $\mathbf{h}(\Sgg)=\log p$, where $p$ is a prime
number. Let $\mu$ be $(\Phi,\s)$-invariant probability measure on
$\Sgg$. Suppose that

\begin{description}

\item[(i)] $\mu$ is ergodic for the action $(\Phi,\s)$;
\item[(ii)] $\mu$ has positive entropy for $\Phi$;
\item[(iii)] the sigma-algebra of the $\s^{(p-1)p}$-invariant sets
coincides $mod\ \mu$ to the sigma-algebra of the $\s$-invariant
sets.
\end{description}

Then, $\mu$ is the Parry measure.}

\cqd

{\rem Given a SC bipermutative c.a. $(\Sgg,\Phi)$, the c.a.
$(K^\Z,\Phi_K)$ obtained from Proposition \ref{IV4} would not be
necessarily bipermutative. For the cases when $(K^\Z,\Phi_K)$ is
bipermutative, we can use Proposition \ref{IV11} to extend for
$(\Sgg,\Phi)$ the results about invariant measures set out by Pivato
\cite{Pivato}.}

\section*{Acknowledgments}

I would like to thank professors A. Maass, S. Mart\'{\i}nez and M.
Pivato, for their discussions and advices on the subject.



\end{document}